\documentclass[12pt,draftcls,onecolumn]{IEEEtran}

\usepackage{amsmath, amsthm, amssymb}
\usepackage{bm}
\usepackage{epsfig}

\newtheorem{thm}{Theorem}[section]

\newtheorem{algo}{Algorithm}[section]

\newtheorem{asmp}{Assumption}[section]

\def\diag{\textrm{diag}}
\def\proof {{Proof.} }
\def\endproof{\hfill $\Box$ \vskip .3cm}
\def\1{\bm 1}
\def\0{\bm 0}

\begin{document}

\title{Limited Attention Allocation in a Stochastic Linear Quadratic System with  Multiplicative Noise}

\author{Xiangyu Cui\thanks{School of Statistics and Management, Shanghai University of Finance and Economics, Shanghai, China}, \quad
	Jianjun Gao\thanks{School of Information Management and Engineering and Research Institute for Interdisciplinary Sciences, Shanghai University of Finance and Economics, Shanghai, China},  \quad
	Lingjie Kong\thanks{School of Statistics and Management, Shanghai University of Finance and Economics, Shanghai, China}}

	
\maketitle

\begin{abstract}
This study addresses limited attention allocation in a stochastic linear quadratic system with multiplicative noise. Our approach enables strategic resource allocation to enhance noise estimation and improve control decisions. We provide analytical optimal control and propose a numerical method for optimal attention allocation. Additionally, we apply our findings to dynamic mean-variance portfolio selection, showing effective resource allocation across time periods and factors, providing valuable insights for investors.

\end{abstract}

\section{Introduction}
The Linear Quadratic (LQ) control problem involving multiplicative noise has been a prominent research area over the past two decades, with significant potential applications in various fields, particularly in mathematical finance. Notably, this problem finds applications in dynamic portfolio selection and financial derivative pricing, as evidenced in studies such as \cite{LiNg}, \cite{zhou2000continuous}, \cite{zhu2004risk}, \cite{primbs2009stochastic}, and \cite{cui2014unified}. Researchers have explored diverse aspects of the LQ control problem with multiplicative noise, investigating formulations with indefinite penalty matrices on the control and state variables (e.g., \cite{lim1999stochastic}, \cite{rami2001solvability}, \cite{costa2007indefinite}). Additionally, scholars have taken keen interest in exploring this problem under various constraints, including cardinality constraint (see \cite{gao2011cardinality}), no shorting constraint (see \cite{li2002dynamic}, \cite{cui2014optimal}), cone constraint (see \cite{czichowsky2013cone}, \cite{cui2017mean}), inequality constraint (see \cite{bemporad2002explicit}), and linear constraint (see \cite{wu2019explicit}).

While many studies on LQ control problem with multiplicative noise assume knowledge of the noise distribution, dynamic portfolio selection practices allow decision makers who employ learning procedures to enhance their estimation of risky assets' returns. This learning process incorporates expert opinions (e.g., \cite{black1991asset}), additional data (e.g., \cite{brennan1998role}, \cite{de2019bayesian}), predictive state variables (e.g., \cite{xia2001learning}), or other information sources (e.g., analysts' recommendations as in \cite{cvitanic2006dynamic}). Bayesian methods are commonly used to update these estimations. These dynamic portfolio selection approaches with learning have considerable relevance to studies on the LQ control problem with parameter uncertainty, which typically focus on additive noise (e.g., \cite{li2008optimal}, \cite{terra2014optimal}, \cite{qian2020dual}). 

\cite{falkinger2008limited} and \cite{gargano2018does} have further demonstrated that during the learning process, decision makers utilize specific attention resources, which are inherently limited. Several works have successfully integrated the concept of limited attention constraints into dynamic portfolio selection (e.g., \cite{peng2006investor}, \cite{andrei2020dynamic}, \cite{zhang2022mean}, \cite{wang2022information}), shedding light on the importance of considering resource constraints in decision-making. However, there is currently a research gap regarding the application of these insights to LQ control problem. Specifically, no existing works have discussed the LQ control problem in the context of limited attention allocation. This unexplored area presents a valuable opportunity to investigate the interplay between limited attention resources and the optimal control policy in stochastic LQ systems. Such analysis could contribute to a deeper understanding of decision-making processes under resource constraints and potentially lead to more effective and realistic control policies.

In this paper, we delve into the limited attention allocation problem within a stochastic linear quadratic control framework with multiplicative noise. Our model combines both an attention allocation policy, responsible for utilizing the limited attention resource to update estimations of the multiplicative noise, and a system control policy, influencing the states of the system. It is essential to note that our studied model is not a mere abstraction of those found in the dynamic portfolio selection literature. While \cite{andrei2020dynamic}, \cite{zhang2022mean}, and \cite{wang2022information} assume that decision makers update risky asset returns using a certain amount of wealth rather than expending attention resources, simplifying the decision problem, our approach explicitly considers the allocation of limited attention resources. Moreover, \cite{peng2006investor} incorporates attention resources for updating risky asset returns, but the study does not address portfolio decision-making or the allocation of attention resources across different periods. Our research significantly contributes to advancing the understanding of the role of limited attention allocation in decision-making for stochastic systems. By analyzing the interplay between limited attention resources, estimation updates, and system control decisions, we aim to provide valuable insights that can enhance decision-making processes in complex constrained environments. 

The paper is organized as follows. In Section 2, we give the basic setting of the problem. In Section 3, we solve the problem, and derive the optimal attention allocation policy and optimal control policy. In Section 4, we explore the practical application of our findings in the context of dynamic portfolio selection. In Section 5, we conclude our paper. 

\section{Basic Setting}
We start our discussion from the following simple LQ control problem,
\begin{align*}
(P_C)~~\min_{\{\bm u_{t+}\}_{t=0}^{T-1}} &~~E_0 \left[ \sum_{t=0}^{T-1} \begin{pmatrix}
\bm u_{t+} & 
x_{t}
\end{pmatrix} \bm B_t \begin{pmatrix}
\bm u_{t+} \\
x_{t}
\end{pmatrix}+ q_T x_T^2 \right]\\
\mbox{s.t.}~~~&~~x_{t+1} = a_{t+1} x_t + \bm b_{t+1}^\prime \bm u_{t+},\\
&~~\bm u_{t+} \in\mathbb{R}^n,\\
&~~\bm b_{t+1} = \bm c + \bm D \bm f_{t+1} + \bm \epsilon_{t+1},\\
&~~\bm f_{t+1} = (\bm I-\bm \varPhi)\bm f_{t} + \bm \eta_{t+1},\\
&~~ x_0 \mbox{ is given},
\end{align*}
where $\bm u_{t+}$ is a $n$-dimensional system control policy, $\bm B_t = \begin{pmatrix}
\bm A_t & \bm p_t\\
\bm p_t^\prime & q_t
\end{pmatrix}$, $\bm A_t\succeq 0$ is a deterministic, symmetric, and semi-positive definite $n\times n$ matrix, $\bm p_t$ is a deterministic $n$-dimensional vector, $q_t\geq 0$ is a deterministic parameter, $q_T>0$ is also a deterministic parameter, $\bm B_t\succeq 0$ is a deterministic, symmetric, and semi-positive definite $(n+1)\times (n+1)$ matrix, $^\prime$ denotes the transpose operator, $a_{t+1}>0$ is a deterministic parameter, $\bm b_{t+1}$ is a $n$-dimensional random vector, $\mathbb{R}^n$ represents a $n$-dimensional Euclidean space, $\bm c$ is a deterministic $n$-dimensional vector, $\bm D$ is a deterministic $n\times k$ matrix, $\bm f_{t+1}$ is a $k$-dimensional random vector ($k<n$), $\bm \epsilon_{t+1}$ follows a $n$-dimensional normal distribution with zero mean and covariance matrix $\bm \Sigma_{\bm \epsilon}$, $\bm I$ is the identity matrix, $\bm \varPhi$ is a deterministic mean-reverting matrix, $\bm \eta_{t+1}$ follows a  $k$-dimensional normal distribution with zero mean and covariance matrix $\bm \Sigma_{\bm \eta}$, and the diagonal elements of $\bm \Sigma_{\bm \eta}$ are denoted as $\sigma_{\eta_1}^2,\cdots, \sigma_{\eta_k}^2$. Different from the classical linear quadratic control problem in \cite{wu2019explicit}, problem $(P_C)$ introduces a factor structure for the random vector $\bm b_{t+1}$ and a mean-reverting structure for the factor $\bm f_{t+1}$ to effectively reduce the dimensionality of the randomness involved. At time $t$, the decision maker can utilize the observed factor $\bm f_t$ to infer the distribution of $\bm b_{t+1}$ and decide the control policy. 

Next, we introduce the limited attention allocation structure into problem $(P_C)$. At time $t$, the decision maker can utilize not only the observed factor $\bm f_t$, but also allocate their attention resources, represented by the scalar parameter $\Lambda_t$, to obtain a noisy signal of $\bm f_{t+1}$,  denoted as $\bm s_{t+}(\bm \lambda_t)$, according to the following relationship, 
\begin{align*}
\bm s_{t+}(\bm \lambda_t) =\bm f_{t+1} + \bm v_{t+1}(\bm \lambda_t),
\end{align*}
where $\bm v_{t+1}(\bm \lambda_t)$ is a white noise with diagonal covariance matrix 
\begin{align*}
\bm \Sigma_{\bm v}(\bm \lambda_t) = \diag \left(\frac{\sigma_{\eta_1}^2}{ e^{ \lambda_{1,t} \theta_1}-1},\dots, \frac{\sigma_{\eta_k}^2}{ e^{ \lambda_{k,t} \theta_k}-1}\right),
\end{align*}
$\bm \lambda_t = (\lambda_{1,t},\cdots,\lambda_{k,t})^{\prime}$ denotes the attention allocation policy, which is the attention resource used for getting the signal, the deterministic parameter $\theta_j$ measures the efficiency of the decision maker’s information processing of the factor $j$.\footnote{The variance setting for noise $\bm v_{t+1}(\bm \lambda_t)$ follows the one proposed in \cite{peng2006investor}. Under the setting, the reduction in the entropy of factor $f_{j,t+1}$ stemming from the knowledge of $s_{j,t+}(\lambda_{j,t})$ is 
	\begin{align*}
	H(f_{j,t+1}) - H(f_{j,t+1}| s_{j,t+}(\lambda_{j,t}))=\frac{1}{2} \lambda_{j,t} \theta_j,
	\end{align*}
	where $H(\cdot)$ is the entropy function, which further implies that the variance of the factor $j$'s posterior distribution is $\sigma_{\eta_j}^2 e^{-\lambda_{j,t} \theta_j}$.}  The larger the $\theta_j$, the higher the efficiency. The attention allocation policy $\bm \lambda_t$ should be nonnegative and no larger than the limited attention resource $\Lambda_t$, i.e., $\bm \lambda_t 
\in \mathbb{R}_+^k,~\bm\lambda_t^T \1 \leq \Lambda_t$. By incorporating the limited attention allocation structure into the linear quadratic control problem, we formulate the linear quadratic control problem with limited attention allocation as follows, 
\begin{align*}
(P_{AC})~~\min_{\{\bm \lambda_t, \bm u_{t+}\}_{t=0}^{T-1}} &~~E_0 \left[ \sum_{t=0}^{T-1} \begin{pmatrix}
\bm u_{t+} & 
x_{t}
\end{pmatrix} \bm B_t \begin{pmatrix}
\bm u_{t+} \\
x_{t}
\end{pmatrix}+q_T x_T^2\right]\\
\mbox{s.t.}~~~&~~x_{t+1} = a_{t+1} x_t + \bm b_{t+1}^\prime \bm u_{t+},\\
&~~\bm u_{t+} \in\mathbb{R}^n,\\
&~~\bm b_{t+1} = \bm c + \bm D \bm f_{t+1} + \bm \epsilon_{t+1}, \\
&~~\bm f_{t+1} = (\bm I-\bm \varPhi)\bm f_{t} + \bm \eta_{t+1}, \\
&~~\Lambda_{t+1} = \Lambda_t - \bm \lambda_t^\prime \bm 1,\\
&~~\bm s_{t+}(\bm \lambda_t) =\bm f_{t+1} + \bm v_{t+1}(\bm \lambda_t),\\
&~~\bm \lambda_t\in \mathbb{R}_+^k,\quad \bm \lambda_t^\prime \bm 1 \leq \Lambda_t,\\
&~~ x_0 \mbox{ and } \Lambda_0 \geq 0 \mbox{ are given}.
\end{align*}
The time $t+$ serves as a dummy time slot used to distinguish between the attention allocation procedure and the system control procedure. It provides a clear separation between these two stages, allowing the decision maker to allocate attention resources at time $t$ and then execute the control policy in the subsequent time slot, $t+$.   All the random vectors are defined in a probability space $(\Omega, \mathcal{F}_T, P)$. The information set at time $t$ is denoted as $\mathcal{F}_t$, which is the $\sigma$-algebra generated by $\{\bm f_0, \bm \epsilon_1,\cdots,\bm \epsilon_t, \bm \eta_1,\cdots,\bm \eta_t, \bm s_{0+}(\bm \lambda_0),\cdots$, $\bm s_{(t-1)+}(\bm \lambda_{t-1})\}$ and the information set after attention allocation at time $t$ is denoted as $\mathcal{F}_{t+}$, which is the $\sigma$-algebra generated by $\{\bm f_0, \bm \epsilon_1,\cdots,\bm \epsilon_t, \bm \eta_1,\cdots,\bm \eta_t$, $\bm s_{0+}(\bm \lambda_0),\cdots$, $\bm s_{t+}(\bm \lambda_{t})\}$.  

\begin{figure}[htb!]
	\centering
	\includegraphics[width=0.85\textwidth]{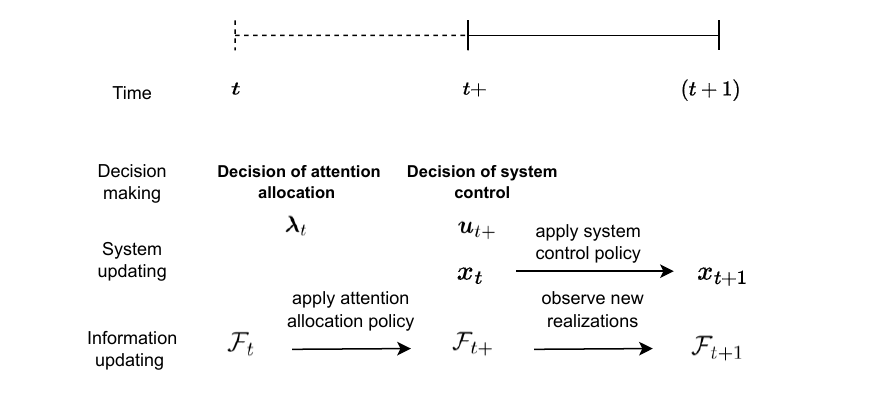}
	\caption{The attention allocation procedure and system control procedure at period $t$}
	\label{fig:period-t}
\end{figure}

To be more intuitive, we display both the attention allocation procedure and system control procedure at period $t$ in Figure \ref{fig:period-t}. During $[t,t+]$, the decision maker chooses the optimal attention allocation policy, generates the signal, and updates the information set and the distribution of the factors. Then, during $[t+,t+1]$, the decision maker chooses the system control policy based on the updated distribution, and updates the information set based on observed new realizations of random vectors. 

Problem $(P_{AC})$ is derived from real-world applications, specifically in the domain of portfolio selection in mathematical finance. In this context, the components of the problem have concrete interpretations: the random vector $\bm b_{t+1}$ represents the excess return of risky assets, the state $x_t$ represents the wealth level at time $t$, and the system control policy represents the portfolio strategy, the factor structure for $\bm b_{t+1}$ follows the well-known factor model, as documented in  \cite{fama1993common} and \cite{herskovic2016common}, and the factors $\bm{f}_{t+1}$ exhibit mean-reverting behavior, a characteristic documented by \cite{garleanu2013dynamic}. The investors can acquire knowledge about the returns of factors using various sources, such as factor valuations (e.g., book-to-market ratio), macroeconomic indicators, cross-sectional characteristics, and machine learning techniques. This process of learning factor returns, known as factor timing, is well-established in the literature (e.g., \cite{Arnott2016}, \cite{Haddad}, \cite{arnott2023factor}, \cite{Favero2021macro}). The economic foundation of the attention allocation structure is built upon this factor timing process.   

At time $t$, the decision maker knows that $\bm f_{t+1}$ and $\bm b_{t+1}$ follow  multivariate normal distributions, 
\begin{align*}
\bm f_{t+1}| \mathcal{F}_t\sim N(\bm \mu_{\bm f,t},\bm \Sigma_{\bm f,t}),\quad \bm b_{t+1}  | \mathcal{F}_{t} \sim N(\bm \mu_{\bm b, t}, \bm \Sigma_{\bm b,t}),
\end{align*}
with
\begin{align*}
&\bm \mu_{\bm f, t}= (\bm I-\bm \varPhi)\bm f_{t},\quad \bm \Sigma_{\bm f, t} = \bm \Sigma_{\bm \eta},\quad \bm \mu_{\bm b, t}= \bm c + \bm D  (\bm I-\bm \varPhi)\bm f_{t},\quad \bm \Sigma_{\bm b, t} = \bm D \bm \Sigma_{\bm \eta}\bm D^{\prime} + \bm \Sigma_{\bm \epsilon}.
\end{align*}
For a given attention allocation strategy $\bm \lambda_t$, the decision maker obtains the signal $\bm s_{t+}(\bm \lambda_t)$ at time $t+$. Then, the decision maker can use the signal to update the conditional distributions of $\bm f_{t+1}$ and $\bm b_{t+1}$. According to Kalman Filtering theory (see Corollary E.3.5 in \cite{Bertsekas}), we have that $\bm f_{t+1}$ and $\bm b_{t+1}$ at time $t+$ are still multivariate normal distributions,
\begin{align*}
\bm f_{t+1}  | \mathcal{F}_{t+} \sim N(\bm \mu_{\bm f, t+}(\bm \lambda_t), \bm \Sigma_{\bm f,t+}(\bm \lambda_t))\quad 
\bm b_{t+1}  | \mathcal{F}_{t+} \sim N(\bm \mu_{\bm b, t+}(\bm \lambda_t), \bm \Sigma_{\bm b,t+}(\bm \lambda_t)),
\end{align*}
where 
\begin{align*}
&\bm \mu_{\bm f, t+} (\bm \lambda_t) = (\bm I-\bm \varPhi)\bm f_{t} +\bm \Sigma_{\bm \eta}(\bm \Sigma_{\bm \eta}+\bm \Sigma_{\bm v}(\bm \lambda_t))^{-1} (\bm s_{t+}(\bm \lambda_t)  - (\bm I-\bm \Phi)\bm f_{t}),\\
&\bm \Sigma_{\bm f,t+}(\bm \lambda_t) = \bm \Sigma_{\bm \eta} - \bm \Sigma_{\bm \eta}(\bm \Sigma_{\bm \eta}+\bm \Sigma_{\bm v}(\bm \lambda_t))^{-1}\bm \Sigma_{\bm \eta},\\
&\bm \mu_{\bm b, t+} (\bm \lambda_t) = \bm c + \bm D(\bm I-\bm \varPhi)\bm f_{t} +\bm D\bm \Sigma_{\bm \eta}(\bm \Sigma_{\bm \eta}+\bm \Sigma_{\bm v}(\bm \lambda_t))^{-1} (\bm s_{t+}(\bm \lambda_t)  - (\bm I-\bm \Phi)\bm f_{t}),\\
&\bm \Sigma_{\bm b,t+}(\bm \lambda_t) = \bm D [\bm \Sigma_{\bm \eta} - \bm \Sigma_{\bm \eta}(\bm \Sigma_{\bm \eta}+\bm \Sigma_{\bm v}(\bm \lambda_t))^{-1}\bm \Sigma_{\bm \eta}]\bm D^\prime + \bm \Sigma_{\bm \epsilon}.
\end{align*}

\section{Optimal Attention Allocation Policy and Optimal Control Policy}
To ensure the problem $(P_{AC})$ has unique optimal control policy, we provide the following assumption.
\begin{asmp}\label{asmp:1}
	For $t=0,1,\cdots,T-1$, the conditional covariance matrix of $\bm b_{t+1}$ at time $t+$, $\bm \Sigma_{\bm b,t+}(\bm \lambda_t)$, is positive definite, i.e., $\bm \Sigma_{\bm b,t+}(\bm \lambda_t)\succ 0$.  
\end{asmp}

The following theorem provides the optimal attention allocation policy and the optimal system control policy for problem $(P_{AC})$. 

\begin{thm}\label{thm:opt}
	Under Assumption \ref{asmp:1}, the optimal control policy and optimal attention allocation policy for problem $(P_{AC})$ at period $t$ are given
	by,
	\begin{align}
	\nonumber &~ \bm u_{t+}^{\star} (x_t, \Lambda_t, \bm f_t, \bm s_{t+}(\bm \lambda_t^\star)) \\
	=&~ -(\bm A_{t}+E_{t+}[h_{t+1}(\Lambda_{t+1}, \bm f_{t+1})\bm b_{t+1}\bm b_{t+1}^\prime] )^{-1} (\bm p_{t}+E_{t+}[h_{t+1}(\Lambda_{t+1}, \bm f_{t+1})\bm b_{t+1}]a_{t+1}) x_t, \\
	\nonumber &~ \bm \lambda_t^{\star}(\Lambda_t, \bm f_t)\\
	\nonumber  =&~ \mathop{\arg\min}_{\substack{\bm \lambda_{t}^\prime\mathbf{1} \leq \Lambda_t, \bm \lambda_{t} \geq \0}}~E_t\big[ E_{t+} [h_{t+1}(\Lambda_{t+1},\bm f_{t+1})] a_{t+1}^2 + q_{t}- (\bm p_{t}+E_{t+}[h_{t+1}(\Lambda_{t+1},\bm f_{t+1})\bm b_{t+1}]a_{t+1})^\prime \\
	&\quad\quad\quad \cdot (\bm A_{t}+E_{t+}[h_{t+1}(\Lambda_{t+1},\bm f_{t+1})\bm b_{t+1}\bm b_{t+1}^\prime] )^{-1} (\bm p_{t}+E_{t+}[h_{t+1}(\Lambda_{t+1},\bm f_{t+1})\bm b_{t+1}]a_{t+1})\big],\label{eq:lambda}
	\end{align}
	where $h_t(\Lambda_t, \bm f_t)$  represents the optimal objective value of the optimization problem (\ref{eq:lambda}) and $h_T(\Lambda_T, \bm f_T)=q_T$, $E_{t+}[\cdot] = E[\cdot|\mathcal{F}_{t+}] $ and $E_{t}[\cdot]=E[\cdot|\mathcal{F}_t]$.
	And the optimal objective function for problem $(P_{AC})$ is given as follows,
	\begin{align*}
	J_0(x_0,\Lambda_0,\bm f_0) = h_0(\Lambda_0, \bm f_0) x_0^2.
	\end{align*}
\end{thm}

\proof
We solve the problem $(P_{AC})$ by dynamic programming. At time $(T-1)+$, the decision maker observes the signal $\bm s_{(T-1)+}(\bm\lambda_{T-1}^{\star})$ and faces the following one-period decision problem,
\begin{align*}
&\min_{\bm u_{(T-1)+}} ~E_{(T-1)+} [\bm u_{(T-1)+}^\prime \bm A_{T-1} \bm u_{(T-1)+} + 2\bm u_{(T-1)+}^\prime \bm p_{T-1} x_{T-1}+ q_{T-1} x_{T-1}^2 +q_T x_T^2]\\
=&\min_{\bm u_{(T-1)+}} ~\bm u_{(T-1)+}^\prime (\bm A_{T-1}+q_{T}E_{(T-1)+}[\bm b_{T}\bm b_{T}^\prime] ) \bm u_{(T-1)+}+ 2\bm u_{(T-1)+}^\prime (\bm p_{T-1}+q_T a_T E_{(T-1)+}[\bm b_{T}]) x_{T-1}  \\
&\quad \quad\quad + q_{T-1} x_{T-1}^2+ q_{T}a_T^2 x_{T-1}^2.
\end{align*}
Under Assumption \ref{asmp:1},  we have $E_{(T-1)+}[\bm b_{T}\bm b_{T}^\prime]\succ 0$. As $q_T> 0$ and $\bm A_{T-1}\succeq 0$, we further have the optimal control policy at period $T-1$ given by,
\begin{align*}
&\bm u_{(T-1)+}^{\star} (x_{T-1},\Lambda_{T-1},\bm f_{T-1}, \bm s_{(T-1)+}(\bm\lambda_{T-1}^{\star})) \\
=& - (\bm A_{T-1}+q_{T}E_{(T-1)+}[\bm b_{T}\bm b_{T}^\prime] )^{-1} (\bm p_{T-1}+q_T a_T E_{(T-1)+}[\bm b_{T}]) x_{T-1},
\end{align*}
where the information set at time $(T-1)+$ contain the signal $\bm s_{(T-1)+}(\bm \lambda_{T-1}^\star)$ and the factor $\bm f_{T-1}$. Substituting the  optimal control policy back, we obtain the cost-to-go function at time $(T-1)+$ as follows,
\begin{align*}
&J_{(T-1)+}(x_{T-1},\Lambda_{T-1},\bm f_{T-1}, \bm s_{(T-1)+}(\bm \lambda_t^{\star}))\\
= &~[q_Ta_T^2 + q_{T-1}- (\bm p_{T-1}+q_T a_T E_{(T-1)+}[\bm b_{T}])^\prime (\bm A_{T-1}+q_{T}E_{(T-1)+}[\bm b_{T}\bm b_{T}^\prime] )^{-1}\\
& \cdot (\bm p_{T-1}+q_T a_T E_{(T-1)+}[\bm b_{T}])] x_{T-1}^2\\
:= &~h_{(T-1)+}(\Lambda_{T-1},\bm \lambda_{T-1}^\star,\bm f_{T-1}) x_{T-1}^2.
\end{align*}
As $\bm \Sigma_{\bm b,(T-1)+}(\bm \lambda_{T-1}^\star)\succ 0$ and $a_T > 0$, which implies $ \begin{pmatrix}
E_{(T-1)+}[\bm b_{T}\bm b_{T}^\prime] & a_T E_{(T-1)+}[\bm b_{T}]\\
a_T E_{(T-1)+}[\bm b_{T}]^\prime & a_T^2
\end{pmatrix}\succ 0$ according to Schur complement. We further have 
\begin{align*}
& \begin{pmatrix}
\bm A_{T-1}+q_{T}E_{(T-1)+}[\bm b_{T}\bm b_{T}^\prime] & \bm p_{T-1}+q_T a_T E_{(T-1)+}[\bm b_{T}]\\
(\bm p_{T-1}+q_T a_T E_{(T-1)+}[\bm b_{T}])^\prime & q_Ta_T^2 + q_{T-1}
\end{pmatrix} \\
=& q_T \begin{pmatrix}
E_{(T-1)+}[\bm b_{T}\bm b_{T}^\prime] & a_T E_{(T-1)+}[\bm b_{T}]\\
a_T E_{(T-1)+}[\bm b_{T}]^\prime & a_T^2
\end{pmatrix} + \bm B_{T-1} \succ 0
\end{align*}
which implies that $h_{(T-1)+}(\bm \lambda_{T-1}^\star,\bm f_{T-1}) >0$ according to Schur complement. Then, at time $(T-1)$, the optimal attention allocation policy is given as,  
\begin{align*}
\bm \lambda_{T-1}^{\star}(\Lambda_{T-1},\bm f_{T-1}) = &\mathop{\arg\min}_{\substack{\bm \lambda_{T-1}^\prime \mathbf{1} \leq \Lambda_{T-1},\bm \lambda_{T-1} \geq \0}}~\{E_{T-1}[J_{(T-1)+}(x_{T-1},\Lambda_{T-1},\bm f_{T-1}, \bm s_{(T-1)+}(\bm \lambda_t))]\}\\
=& \mathop{\arg\min}_{\substack{\bm \lambda_{T-1}^\prime \mathbf{1} \leq \Lambda_{T-1},\bm \lambda_{T-1} \geq \0}}~~\{E_{T-1}[h_{(T-1)+}(\Lambda_{T-1},\bm \lambda_{T-1},\bm f_{T-1})]\}.
\end{align*}
And the cost-to-go function at time $(T-1)$ becomes
\begin{align*}
J_{T-1}(x_{T-1},\Lambda_{T-1},\bm f_{T-1} )=h_{T-1}(\Lambda_{T-1},\bm f_{T-1}) x_{T-1}^2.
\end{align*}

We assume that the cost-to-go function at time $t+1$ is $J_{t+1}(x_{t+1},\Lambda_{t+1},\bm f_{t+1} )=h_{t+1}(\Lambda_{t+1},\bm f_{t+1}) x_{t+1}^2$ and $h_{t+1}(\Lambda_{t+1},\bm f_{t+1})>0$. Then, we derive the optimal control policy and attention allocation policy for period $t$. At time $t+$,  the investors decision maker observes the signal $\bm s_{t+}(\bm\lambda_t^{\star})$ and face the following one-period portfolio selection problem according to Bellman's principal of optimality,
\begin{align*}
&\min_{\bm u_{t+}} ~E_{t+} [\bm u_{t+}^\prime \bm A_{t} \bm u_{t+} + 2\bm u_{t+}^\prime \bm p_{t} x_{t} + q_{t} x_{t}^2+J_{t+1}(x_{t+1},\Lambda_{t+1},\bm f_{t+1})]\\
=&\min_{\bm u_{t+}} ~\bm u_{t+}^\prime (\bm A_{t}+E_{t+}[h_{t+1}(\Lambda_{t+1},\bm f_{t+1})\bm b_{t+1}\bm b_{t+1}^\prime] ) \bm u_{t+}\\
&\quad \quad ~~ + 2\bm u_{t+}^\prime (\bm p_{t}+E_{t+}[h_{t+1}(\Lambda_{t+1},\bm f_{t+1})\bm b_{t+1}]a_{t+1}) x_{t} + q_{t} x_{t}^2 + E_{t+} [h_{t+1}(\Lambda_{t+1},\bm f_{t+1})] a_{t+1}^2 x_{t}^2.
\end{align*}
We define the new probability measure $Q_{t+}(\bm \lambda_{t}^\star)$ after observing the signal $\bm s_{t+}(\bm \lambda_{t}^\star)$, which is equivalent to the objective probability measure $P$, as follows,
\begin{align*}
\frac{dQ_{t+}(\bm \lambda_{t}^\star)}{dP} = \frac{h_{t+1}(\Lambda_{t+1},\bm f_{t+1})}{E_{t+}[h_{t+1}(\Lambda_{t+1},\bm f_{t+1})]}.
\end{align*}
As the new probability measure $Q_{t+}(\bm \lambda_{t}^\star)$ is equivalent to $P$, the covariance matrix of $\bm b_{t+1}$ under $Q_{t+}(\bm \lambda_{t}^\star)$ is also positive definite, which implies that $\bm A_{t}+E_{t+}[h_{t+1}(\Lambda_{t+1},\bm f_{t+1})\bm b_{t+1}\bm b_{t+1}^\prime]$ is positive definite. Thus, the optimal control policy at period $t$ is 
\begin{align*}
&\bm u_{t+}^{\star}(x_t, \Lambda_t, \bm f_t, \bm s_{t+}(\bm \lambda_t^\star)) \\
= &-(\bm A_{t}+E_{t+}[h_{t+1}(\Lambda_{t+1},\bm f_{t+1})\bm b_{t+1}\bm b_{t+1}^\prime] )^{-1} (\bm p_{t}+E_{t+}[h_{t+1}(\Lambda_{t+1},\bm f_{t+1})\bm b_{t+1}]a_{t+1}) x_t,
\end{align*}
and the cost-to-go function at time $t+$ becomes
\begin{align*}
& J_{t+}(x_t, \Lambda_t,\bm f_{t}, \bm s_{t+}(\bm \lambda_t^{\star}))\\
=& \{E_{t+} [h_{t+1}(\Lambda_{t+1},\bm f_{t+1})] a_{t+1}^2 + q_{t}- (\bm p_{t}+E_{t+}^{Q_{t+}(\bm \lambda_t^\star)}[\bm b_{t+1}] E_{t+} [h_{t+1}(\Lambda_{t+1},\bm f_{t+1})] a_{t+1})^\prime \\
& \cdot (\bm A_{t}+E_{t+}^{Q_{t+}(\bm \lambda_t^\star)}[\bm b_{t+1}\bm b_{t+1}^\prime] E_{t+} [h_{t+1}(\Lambda_{t+1},\bm f_{t+1})] )^{-1}\\
& \cdot (\bm p_{t}+E_{t+}^{Q_{t+}(\bm \lambda_t^\star)}[\bm b_{t+1}]E_{t+} [h_{t+1}(\Lambda_{t+1},\bm f_{t+1})]a_{t+1})\} x_t^2\\
:=& h_{t+}(\Lambda_t,\bm \lambda_t^{\star},\bm f_t)  x_{t}^2,
\end{align*}
where $E_{t+}^{Q_{t+}(\bm \lambda_t^\star)}[\cdot]$ denotes the conditional expectation under probability $Q_{t+}(\bm \lambda_t^\star)$. Following the same argument for $T-1$ period, we have 
\begin{align*}
& \begin{pmatrix}
\bm A_{t}+E_{t+}^{Q_{t+}(\bm \lambda_t^\star)}[\bm b_{t+1}\bm b_{t+1}^\prime] E_{t+} [h_{t+1}(\Lambda_{t+1},\bm f_{t+1})]  & \bm p_{t}+E_{t+}^{Q_{t+}(\bm \lambda_t^\star)}[\bm b_{t+1}]E_{t+} [h_{t+1}(\Lambda_{t+1},\bm f_{t+1})]a_{t+1}\\
(\bm p_{t}+E_{t+}^{Q_{t+}(\bm \lambda_t^\star)}[\bm b_{t+1}]E_{t+} [h_{t+1}(\Lambda_{t+1},\bm f_{t+1})]a_{t+1})^\prime &  q_{t} + E_{t+} [h_{t+1}(\Lambda_{t+1},\bm f_{t+1})]a_{t+1}^2
\end{pmatrix} \\
=& E_{t+} [h_{t+1}(\Lambda_{t+1}, \bm f_{t+1})] \begin{pmatrix}
E_{t+}^{Q_{t+}(\bm \lambda_t^\star)}[\bm b_{t+1}\bm b_{t+1}^\prime] & a_{t+1} E_{t+}^{Q_{t+}(\bm \lambda_t^\star)}[\bm b_{t+1}]\\
a_{t+1}E_{t+}^{Q_{t+}(\bm \lambda_t^\star)}[\bm b_{t+1}^\prime] & a_{t+1}^2
\end{pmatrix} + \bm B_{t} \succ 0
\end{align*}
which implies that $h_{t+}(\Lambda_{t},\bm \lambda_{t}^\star,\bm f_{t}) >0$ according to Schur complement.  Then, at time $t$, the optimal attention allocation strategy  is given as,
\begin{align*}
\bm \lambda_{t}^{\star}(\Lambda_t,\bm f_t) = &\mathop{\arg\min}_{\substack{\bm \lambda_{t}^\prime \mathbf{1} \leq \Lambda_{t},
		\bm \lambda_{t} \geq \0}}~\{E_{t}[J_{t+}(x_t, \Lambda_t,\bm f_{t}, \bm s_{t+}(\bm \lambda_t))]\}=\mathop{\arg\min}_{\substack{\bm \lambda_{t}^\prime \mathbf{1} \leq \Lambda_{t}, 
		\bm \lambda_{t} \geq \0}}~\{E_{t}[h_{t+}(\Lambda_t, \bm \lambda_{t},\bm f_t)]\}.
\end{align*}
And the cost-to-go function at time $t$ becomes
\begin{align*}
J_{t}(x_t,\Lambda_t,\bm f_{t} )=h_{t}(\Lambda_t,\bm f_t) x_{t}^2.
\end{align*}

Therefore, we have prove our main results. Furthermore, the optimal objective of problem is $J_0(x_0,\Lambda_0,\bm f_0) = h_0(\Lambda_0,\bm f_0) x_0^2$.
\endproof

According to the theorem, several significant findings emerge. First, the optimal control policy is a linear function of the current state $x_t$, implying that the decision maker's system control policy is determined based on the present state. On the other hand, the optimal attention allocation policy is independent of the current state $x_t$ but is influenced by the attention resource $\Lambda_t$ and the factors $\bm{f}_t$. The decision maker allocates attention resources based on the available attention capacity and the current factors' information. Second, the optimization problem used to determine the optimal attention allocation policy is not convex, indicating that finding the optimal allocation requires numerical techniques rather than a straightforward analytical solution. Third, both the optimal control policy and the optimal attention allocation policy heavily rely on the random variable $h_t(\Lambda_t, \bm{f}_t)$. This variable represents the potential future control opportunity, reflecting the expected impact of future decisions on current decision-making process. Fourth, when no attention allocation decisions are made, the optimal control policy reduces to the policy proposed in Proposition 1 of \cite{wu2019explicit}. This comparison highlights the significance of incorporating attention allocation decisions to linear quadratic control problem.

Next, we propose the following algorithm to derive $h_t(\Lambda_t,\bm f_t)$ and the optimal attention allocation policy $\bm \lambda_t^{\star}(\Lambda_t, \bm f_t)$.
\begin{algo}\label{algo:1}
\begin{enumerate}
	\item Choose the constant function, $h_T(\Lambda_T,\bm f_T) = q_T$ and set $t=T-1$. Simulate $L$ realizations of the independent multivariate standard normal random vectors, $\bm \varepsilon_1^{(\ell_1)}$ ($k$-dimensional), $\bm \varepsilon_2^{(\ell_2)}$ ($k$-dimensional), $\bm \varepsilon_3^{(\ell_2)}$ ($n$-dimensional).  
	\item Set up a finite grid consisting of $N$ values of $\Lambda_t$ and $M$ realizations of $\bm f_t$.  
	\item For each realization $\bm f_{t}^{(i)}$ and possible $\Lambda_t^{(j)}$, compute $\bm \lambda_{t}^\star(\Lambda_t^{(j)}, \bm f_{t}^{(i)})$ by solving the optimization problem numerically,  
	\begin{align*}
	&\min_{\bm \lambda_{t}^\prime \1\leq \Lambda_t^{(j)},~\bm \lambda_t \geq \0 } \frac{1}{L} \sum_{\ell_1=1}^{L} \bigg\{q_t+ \frac{a_{t+1}^2}{L} \sum_{\ell_2=1}^L h_{t+1} \big(\Lambda_{t+1}^{(j)}(\bm \lambda_{t}),\bm f_{t+1}^{(\ell_1,\ell_2,i)}\big) \\
	& \quad - \bigg(\bm p_t + \frac{a_{t+1}}{L} \sum_{\ell_2=1}^L h_{t+1} (\Lambda_{t+1}^{(j)}(\bm \lambda_{t}),\bm f_{t+1}^{(\ell_1,\ell_2,i)}) \bm b_{t+1}^{(\ell_1,\ell_2,i)}\bigg)^\prime \\
	&\quad \cdot\bigg(\bm A_t + \frac{1}{L} \sum_{\ell_2=1}^L h_{t+1} \big(\Lambda_{t+1}^{(j)}(\bm \lambda_{t}),\bm f_{t+1}^{(\ell_1,\ell_2,i)}\big) \bm b_{t+1}^{(\ell_1,\ell_2,i)} (\bm b_{t+1}^{(\ell_1,\ell_2,i)})^\prime \bigg)^{-1}\\
	&\quad \cdot\bigg(\bm p_t + \frac{a_{t+1}}{L} \sum_{\ell_2=1}^L h_{t+1} \big(\Lambda_{t+1}^{(j)}(\bm \lambda_{t}),\bm f_{t+1}^{(\ell_1,\ell_2,i)}\big) \bm b_{t+1}^{(\ell_1,\ell_2,i)}\bigg)\bigg\},\\
	&:= \min_{\bm \lambda_{t}^\prime \1\leq \Lambda_t^{(j)},~\bm \lambda_t \geq \0 } \frac{1}{L} \sum_{\ell_1=1}^{L} g(\bm \lambda_t; \Lambda_t^{(j)}, \bm f_t^{(i)}, \bm s_{t+}(\bm \lambda_t)^{(\ell_1,i)}),
	\end{align*}
	where 
	\begin{align*}
	&\Lambda_{t+1}^{(j)}(\bm \lambda_{t})=\Lambda_t^{(j)}-\bm \lambda_{t}^\prime \1,\\
	&\bm f_{t+1}^{(\ell_1,\ell_2,i)} = (\bm I-\bm \varPhi)\bm f_{t}^{(i)} + \bm \Sigma_{\bm \eta}(\bm \Sigma_{\bm \eta}+\bm \Sigma_{\bm v}(\bm \lambda_t))^{-1} (\bm s_{t+}(\bm \lambda_t)^{(\ell_1,i)}  - (\bm I-\bm \Phi)\bm f_{t}^{(i)}) \\
	&\quad \quad \quad \quad + [\bm \Sigma_{\bm \eta} - \bm \Sigma_{\bm \eta}(\bm \Sigma_{\bm \eta}+\bm \Sigma_{\bm v}(\bm \lambda_t))^{-1}\bm \Sigma_{\bm \eta}]^{1/2} \bm \varepsilon_2^{(\ell_2)},\\
	& \bm b_{t+1}^{(\ell_1,\ell_2,i)} = \bm c +\bm D \bm f_{t+1}^{(\ell_1,\ell_2,i)} + \bm \Sigma_{\bm \epsilon}^{1/2} \bm \varepsilon_3^{(\ell_2)},\\
	&  \bm s_{t+}(\bm \lambda_t)^{(\ell_1,i)} = (\bm I-\bm \varPhi)\bm f_{t}^{(i)} + ( \bm\Sigma_{\bm \eta} + \bm \Sigma_{\bm \nu}(\bm \lambda_t)) \bm \varepsilon_1^{(\ell_1)}.
	\end{align*}	
	Denote the optimal value as $h_t(\Lambda_t^{(j)},\bm f_t^{(i)})$ and the optimizer as $\bm \lambda_t^\star(\Lambda_t^{(j)},\bm f_t^{(i)})$.   
	\item Approximate the functions $h_t(\Lambda_t, \bm f_{t})$ and $\bm \lambda_t^\star(\Lambda_t,\bm f_t)$ by linear interpolation. Set $t=t-1$. 
	\item When $t=0$, stop the computation. When $t>0$, go back to step 2.
\end{enumerate}
\label{alg:1}
\end{algo}

In step 3 of Algorithm \ref{alg:1}, the objective of the optimization problem is not a simple summation over all simulated random variables. For a given realization of $\bm \varepsilon_1^{(\ell_1)}$, we can compute the gradient of function $g(\bm \lambda_t; \Lambda_t^{(j)}, \bm f_t^{(i)}, \bm s_{t+}(\bm \lambda_t)^{(\ell_1,i)})$ with respect to $\bm \lambda_t$ and apply the stochastic gradient decent method to solve the optimization problem. Please note that the complexity of Algorithm \ref{alg:1} does not grow exponentially with respect to period $T$.  

\section{Application in Dynamic Mean-Variance Portfolio Selection}
In this section, we explore an application of our findings in the context of dynamic portfolio selection and consider a discrete-time mean-variance model with limited attention allocation. In this model, the investor's objective is to minimize the variance of the terminal wealth $x_T$ while ensuring that the expected terminal wealth achieves a target value $d \geq \rho_0 x_0$, where $a_{t+1}>1$ is the total return of a risk-free investment at period $t$ and $\rho_t=\prod_{s=t}^{T-1}a_{s+1}$, i.e.,
\begin{align*}
\min_{\{\bm \lambda_t,\bm u_{t+}\}_{t=0}^{T-1}} & \mbox{Var}_0 (x_T)\\ \mbox{s.t.}~~~& E_0[x_T] = d,\\
& (x_t, \bm \lambda_t, \Lambda_t, \bm u_{t+}) \mbox{ obey constraints of } (P_{AC}),
\end{align*}
where $ \mbox{Var}_0 (x_T)$ denotes the variance of the terminal wealth. Following \cite{li2002dynamic}, we introduce the Lagrangian multiplier $2\mu$ for the constraint of expected terminal wealth, define $y_t = x_t - (d-\mu) \rho_t^{-1}$, and obtain the auxiliary problem as follows, 
\begin{align*}
 (\mathcal{A}_\mu)~~\min_{\{\bm \lambda_t,\bm u_{t+}\}_{t=0}^{T-1}} &~ E_0[y_T^2] - \mu^2\\
\mbox{s.t.}~~~& (y_t, \bm \lambda_t, \Lambda_t, \bm u_{t+}) \mbox{ obey constraints of } (P_{AC}).
\end{align*} 
Then, based on strong duality, we can solve the dynamic mean-variance problem $(MV)$ by optimizing the objective function of $(\mathcal{A}_\mu)$ over the set $\{\mu|\mu\in \mathbb{R}\}$. 

\begin{thm}
	Under the assumption that $E_{(T-1)+}[\bm b_{T}]\neq \0$, the optimal portfolio policy  and optimal attention allocation policy are given by
	\begin{align}
	\nonumber &~ \bm u_{t+}^{\star} (x_t, \Lambda_t, \bm f_t, \bm s_{t+}(\bm \lambda_t^\star)) \\
	=&~ - E_{t+}[h_{t+1}(\Lambda_{t+1}, \bm f_{t+1})\bm b_{t+1}\bm b_{t+1}^\prime] ^{-1} E_{t+}[h_{t+1}(\Lambda_{t+1}, \bm f_{t+1})\bm b_{t+1}] a_{t+1} (x_t-\rho_t(d-\mu^\star)), \\
	\nonumber &~ \bm \lambda_t^{\star}(\Lambda_t, \bm f_t)\\
	\nonumber  =&~ \mathop{\arg\min}_{\substack{\bm \lambda_{t}^\prime \mathbf{1} \leq \Lambda_t, \bm \lambda_{t} \geq \0}}~  a_{t+1}^2 E_t\big[ E_{t+} [h_{t+1}(\Lambda_{t+1},\bm f_{t+1})] - E_{t+}[h_{t+1}(\Lambda_{t+1},\bm f_{t+1})\bm b_{t+1}^\prime] \\
	&\quad\quad\quad \cdot E_{t+}[h_{t+1}(\Lambda_{t+1},\bm f_{t+1})\bm b_{t+1}\bm b_{t+1}^\prime]^{-1} E_{t+}[h_{t+1}(\Lambda_{t+1},\bm f_{t+1})\bm b_{t+1}]\big],\label{eq:lambda-2}
	\end{align}
	where $h_t(\Lambda_t, \bm f_t)$ represents the optimal objective value of the optimization problem (\ref{eq:lambda-2}) and $h_T(\Lambda_T, \bm f_T)=1$, and the optimal multiplier $\mu^\star =  \frac{h_{0}(\Lambda_0,\bm f_0) (\rho_0x_0 - d)}{\rho_0^2-h_{0}(\Lambda_0,\bm f_0)}$. 
	The efficient frontier is given as follows,
	\begin{align*}
	\mbox{Var}_0 (x_T) = \frac{h_{0}(\Lambda_0,\bm f_0)  (E_0[x_T] -\rho_0x_0 )^2}{\rho_0^2-h_{0}(\Lambda_0,\bm f_0)  },\quad E_0[x_T] \geq \rho_0 x_0.
	\end{align*}
\end{thm}

\proof
The auxiliary problem is a particular linear quadratic control system with limited attention allocation. With the help of the results in Theorem \ref{thm:opt}, we only need to compute the optimal Lagrangian multiplier $\mu^{\star}$ according to
\begin{align*}
\mu^{\star} &= \mathop{\arg\max} \limits_{\mu \in \mathbb{R}}~~ \{h_0(\Lambda_0, \bm f_0) (x_0 - \rho_0^{-1}(d-\mu))^2 - \mu^2\}.
\end{align*}
As $E_{(T-1)+}[\bm b_{T}\bm b_{T}^\prime]^{-1}\succ 0$ and $E_{(T-1)+}[\bm b_{T}]\neq \0$, we have 
$h_{T-1}(\Lambda_{T-1},\bm f_{T-1}) < a_{T-1}^2$. With the help of $E_{t+}[h_{t+1}(\Lambda_{t+1},\bm f_{t+1})\bm b_{t+1}\bm b_{t+1}^\prime]^{-1}\succ 0$, we further have $h_{t}(\Lambda_{t},\bm f_{t}) < \rho_t^2$. Thus, the optimal Lagrangian multiplier $\mu^{\star}$ is
$\mu^{\star}= \frac{h_{0}(\Lambda_0,\bm f_0) (\rho_0x_0 - d)}{\rho_0^2-h_{0}(\Lambda_0,\bm f_0)}$.
Moreover, the variance of the terminal wealth level achieved by optimal portfolio policy is
\begin{align*}
\mbox{Var}_0(x_T) & = h_{0}(\Lambda_0,\bm f_0) \left(x_0 -\rho_0^{-1}d+\rho_0^{-1}\mu^{\star}\right)^2 - (\mu^{\star})^2=\frac{h_{0}(\Lambda_0,\bm f_0)  (\rho_0x_0 - d)^2}{\rho_0^2-h_{0}(\Lambda_0,\bm f_0)  },
\end{align*}
which provides the mean-variance efficient frontier by changing $d$ into $E_0[x_T]$.
\endproof

In our empirical analysis, we utilize data from the U.S. market obtained from Kenneth R. French's data library website. The dataset comprises monthly returns of the Fama-French three factors, denoted as $\bm{f}_t=(MKT_t, SMB_t, HML_t)^{\prime}$, and the monthly returns of U.S. five industry portfolios, including Cnsmr (C), Manuf (Ma), HiTec (Hi), Hlth (Hl), and Other (O), represented as $\bm{b}_{t}=(b_{C,t}, b_{Ma,t},b_{Hi,t}, b_{Hl,t},b_{O,t})^{\prime}$. We set $a_{t+1}$ equal to the average total return of a risk-free investment, $1.0036$.


In our analysis, we make the assumption that the efficiency parameters of the decision maker's information processing of the Fama-French three factors are all equal to 0.69, i.e., $\theta_{MKT}=\theta_{SMB}=\theta_{HML} = 0.69$.\footnote{When the decision maker uses 1 unit attention resource to learn a factor, the conditional variance of the factor under $\mathcal{F}_{t+}$ is half of the conditional variance of the factor under $\mathcal{F}_t$.} Under the parameter setting, we consider a three-month mean-variance model with limited attention allocation. The initial attention resource is $\Lambda_0=3$. As the optimal portfolio policy is analytical, we mainly focus on investigating the optimal attention allocation policy. Applying Algorithm \ref{algo:1}, we compute the optimal attention allocation policy, $\bm \lambda_{t}^\star(\Lambda_t,\bm f_t)$ for $t=0,1,2$, with $13\times7 \times 7\times 7$ grid points. 

\begin{figure}[htb!]
\centering
\includegraphics[width=0.85\textwidth] {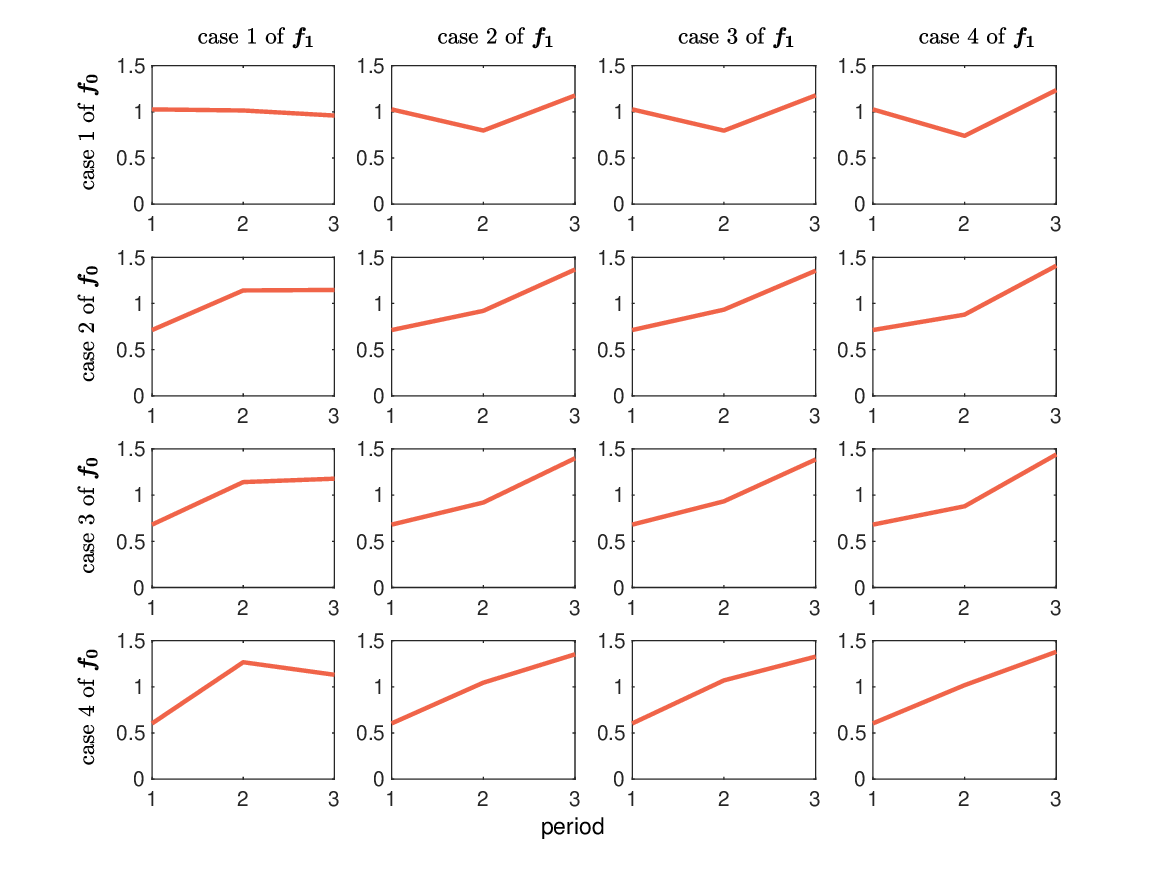}   
\caption{The total optimal attention allocation in different cases of factors and different periods}
\label{fig2}
\end{figure}

We first discuss the optimal attention allocation among different periods. Figure \ref{fig2} shows the total optimal attention allocation $\1^\prime \bm \lambda_{t}^\star$ in different cases of factors and different periods. Case 1 to case 4 correspond to four different settings for the factors $\bm{f}_t$:  $(0.005,0.002,0.003)^\prime$, $(-0.127,0.002,0.003)^\prime$, $(-0.127,-0.091,0.003)^\prime$, $(-0.127,-0.091,-0.081)^\prime$.\footnote{$(0.005,0.002,0.003)^\prime$ is the mean of $\bm f_t$, and $(-0.127,-0.091,-0.081)^\prime$ is the mean minus three times the standard  deviation of $\bm f_t$ for each factor. The cases with positive extreme values have the similar pattern, which are not reported in the paper.} Based on the total optimal attention allocation for period 1 in the four cases of each column, we can observe an interesting pattern. When the values of the factors are extreme values, particularly in Case 4 where the market state is three standard deviations below the mean for each factor, the mean-reverting effect of the factors may generate relatively accurate estimations of the future factor. As a result, the investors tend to allocate less attention resource to learn the factors. The total optimal attention allocation for period 2 in the four cases of each row confirm the finding.  In the last period, the decision maker faces a critical choice as there is no future period to allocate attention to. As a result, the decision maker may opt to use all the remaining attention resource available to make the most informed investment decisions in the current period, which leads to lower total optimal attention allocation for the left and top cases compared to the right and bottom cases.

\begin{figure}[htb!]
	\centering
	\includegraphics[width=0.85\textwidth] {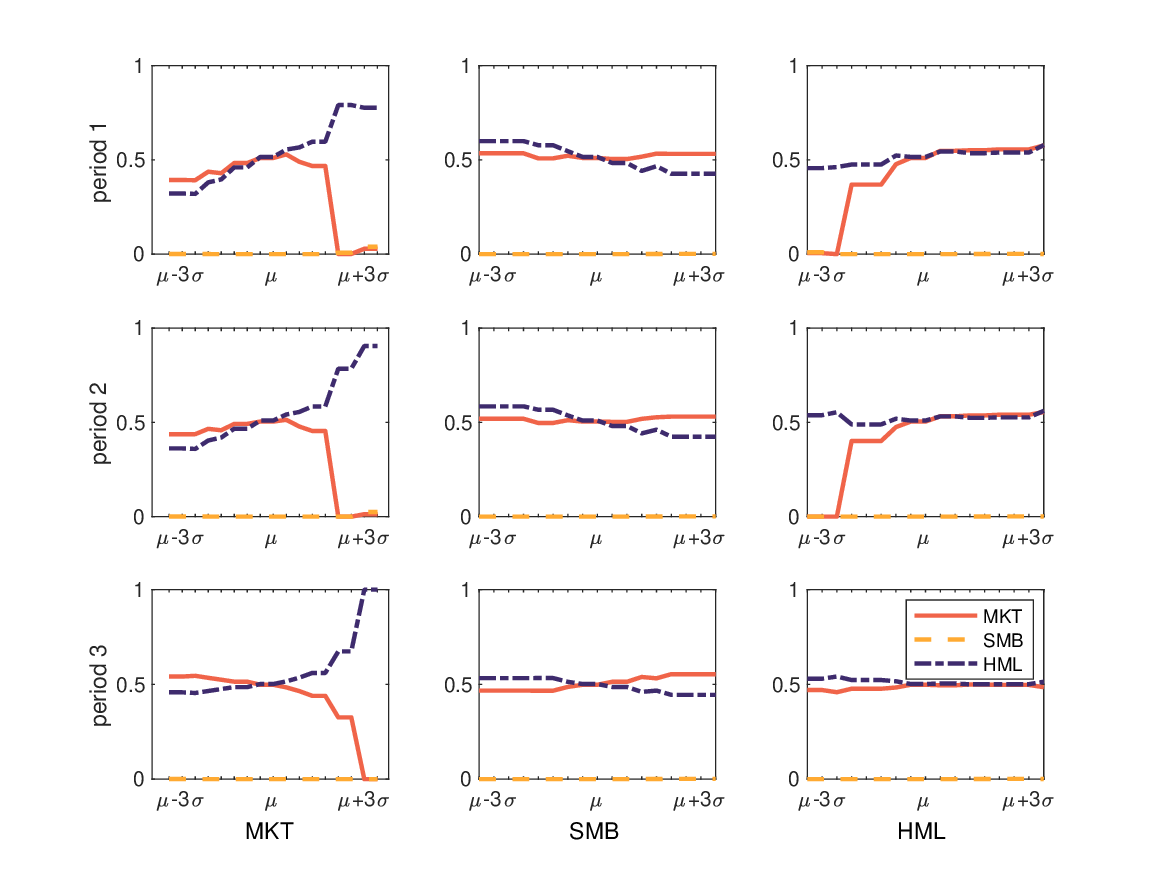}   
	\caption{The optimal attention allocation among different factors} 
	\label{fig3}
\end{figure}

Next, we examine optimal attention allocation among factors in Figure \ref{fig3} when $\bm f_0$ and $\bm f_1$ are both from Case 1. Each curve in the plot represents the allocation in different factors over three periods. Investors prioritize attention on MKT and HML factors. Although the first and second periods face the same factors, their allocations differ due to varying resource availability and future opportunities. In the last period, with no future opportunities, investors solely focus on the single-period problem. The figure reveals how investors strategically allocate attention, considering market conditions, resource availability, and future opportunities, influencing their decisions over multiple time periods.

Finally, during the out-of-sample period from Jan. 2008 to Dec. 2017, we apply a rolling window scheme to estimate parameters and conduct three months mean-variance portfolio selection with limited attention allocation. This results in 118 terminal wealth values. Without attention allocation, the average terminal wealth is $1.0648$ (standard deviation: $0.2177$, Sharpe ratio: $0.2477$). With attention allocation, the average terminal wealth improves to $1.0848$ (standard deviation: $0.2694$, Sharpe ratio: $0.2744$), a 10.8\% boost in the Sharpe ratio. The strategy's superior risk-adjusted returns highlight its ability to make more informed investment decisions.

\section{Conclusion}

This paper introduces limited attention allocation into the classical linear quadratic control problem. The model allows strategic resource allocation, improving randomness estimations. We derive explicit optimal control and a numerical method for attention allocation. Empirical analysis extends the model to dynamic portfolio selection, showcasing effective resource allocation over time and factors. Attention allocation significantly enhances decision-making and risk-adjusted performance. The study provides valuable insights into the importance of limited attention allocation in decision-making processes.

Our model can be extended for real-world application. We can replace the abstract white noise setting with concrete signal generation based on factor timing literature, enhancing applicability. The decision maker can also allocate resources using optimal weights from pre-given prediction models, providing flexibility beyond factor attention. These extensions improve the model's practical relevance and empower investors to make informed decisions in diverse scenarios.

\bibliographystyle{ieeetr}
\bibliography{ref}

\begin{thebibliography}{10}

\bibitem{LiNg}
D.~Li and W.~L. Ng, ``Optimal dynamic portfolio selection: Multiperiod
  mean-variance formulation,'' {\em Mathematical Finance}, vol.~10, no.~3,
  pp.~387--406, 2000.

\bibitem{zhou2000continuous}
X.~Y. Zhou and D.~Li, ``Continuous-time mean-variance portfolio selection: A
  stochastic lq framework,'' {\em Applied Mathematics and Optimization},
  vol.~42, pp.~19--33, 2000.

\bibitem{zhu2004risk}
S.~S. Zhu, D.~Li, and S.~Y. Wang, ``Risk control over bankruptcy in dynamic
  portfolio selection: A generalized mean-variance formulation,'' {\em IEEE
  Transactions on Automatic Control}, vol.~49, no.~3, pp.~447--457, 2004.

\bibitem{primbs2009stochastic}
J.~A. Primbs and C.~H. Sung, ``Stochastic receding horizon control of
  constrained linear systems with state and control multiplicative noise,''
  {\em IEEE transactions on Automatic Control}, vol.~54, no.~2, pp.~221--230,
  2009.

\bibitem{cui2014unified}
X.~Cui, X.~Li, and D.~Li, ``Unified framework of mean-field formulations for
  optimal multi-period mean-variance portfolio selection,'' {\em IEEE
  Transactions on Automatic Control}, vol.~59, no.~7, pp.~1833--1844, 2014.

\bibitem{lim1999stochastic}
A.~E. Lim and X.~Y. Zhou, ``Stochastic optimal lqr control with integral
  quadratic constraints and indefinite control weights,'' {\em IEEE
  Transactions on Automatic Control}, vol.~44, no.~7, pp.~1359--1369, 1999.

\bibitem{rami2001solvability}
M.~A. Rami, X.~Chen, J.~B. Moore, and X.~Y. Zhou, ``Solvability and asymptotic
  behavior of generalized riccati equations arising in indefinite stochastic lq
  controls,'' {\em IEEE Transactions on Automatic Control}, vol.~46, no.~3,
  pp.~428--440, 2001.

\bibitem{costa2007indefinite}
O.~L. Costa and W.~L. de~Paulo, ``Indefinite quadratic with linear costs
  optimal control of markov jump with multiplicative noise systems,'' {\em
  Automatica}, vol.~43, no.~4, pp.~587--597, 2007.

\bibitem{gao2011cardinality}
J.~Gao and D.~Li, ``Cardinality constrained linear-quadratic optimal control,''
  {\em IEEE Transactions on Automatic Control}, vol.~56, no.~8, pp.~1936--1941,
  2011.

\bibitem{li2002dynamic}
X.~Li, X.~Y. Zhou, and A.~E. Lim, ``Dynamic mean-variance portfolio selection
  with no-shorting constraints,'' {\em SIAM Journal on Control and
  Optimization}, vol.~40, no.~5, pp.~1540--1555, 2002.

\bibitem{cui2014optimal}
X.~Cui, J.~Gao, X.~Li, and D.~Li, ``Optimal multi-period mean--variance policy
  under no-shorting constraint,'' {\em European Journal of Operational
  Research}, vol.~234, no.~2, pp.~459--468, 2014.

\bibitem{czichowsky2013cone}
C.~Czichowsky and M.~Schweizer, ``Cone-constrained continuous-time markowitz
  problems,'' {\em The Annals of Applied Probability}, vol.~23, no.~2,
  pp.~764--810, 2013.

\bibitem{cui2017mean}
X.~Cui, D.~Li, and X.~Li, ``Mean-variance policy for discrete-time
  cone-constrained markets: Time consistency in efficiency and the
  minimum-variance signed supermartingale measure,'' {\em Mathematical
  Finance}, vol.~27, no.~2, pp.~471--504, 2017.

\bibitem{bemporad2002explicit}
A.~Bemporad, M.~Morari, V.~Dua, and E.~N. Pistikopoulos, ``The explicit linear
  quadratic regulator for constrained systems,'' {\em Automatica}, vol.~38,
  no.~1, pp.~3--20, 2002.

\bibitem{wu2019explicit}
W.~Wu, J.~Gao, D.~Li, and Y.~Shi, ``Explicit solution for constrained
  scalar-state stochastic linear-quadratic control with multiplicative noise,''
  {\em IEEE Transactions on Automatic Control}, vol.~64, no.~5, pp.~1999--2012,
  2019.

\bibitem{black1991asset}
F.~Black and R.~B. Litterman, ``Asset allocation: Combining investor views with
  market equilibrium,'' {\em The Journal of Fixed Income}, vol.~1, no.~2,
  pp.~7--18, 1991.

\bibitem{brennan1998role}
M.~J. Brennan, ``The role of learning in dynamic portfolio decisions,'' {\em
  Review of Finance}, vol.~1, no.~3, pp.~295--306, 1998.

\bibitem{de2019bayesian}
C.~De~Franco, J.~Nicolle, and H.~Pham, ``Bayesian learning for the markowitz
  portfolio selection problem,'' {\em International Journal of Theoretical and
  Applied Finance}, vol.~22, no.~07, p.~1950037, 2019.

\bibitem{xia2001learning}
Y.~Xia, ``Learning about predictability: The effects of parameter uncertainty
  on dynamic asset allocation,'' {\em The Journal of Finance}, vol.~56, no.~1,
  pp.~205--246, 2001.

\bibitem{cvitanic2006dynamic}
J.~Cvitani{\'c}, A.~Lazrak, L.~Martellini, and F.~Zapatero, ``Dynamic portfolio
  choice with parameter uncertainty and the economic value of analysts’
  recommendations,'' {\em The Review of Financial Studies}, vol.~19, no.~4,
  pp.~1113--1156, 2006.

\bibitem{li2008optimal}
D.~Li, F.~Qian, and P.~Fu, ``Optimal nominal dual control for discrete-time
  linear-quadratic gaussian problems with unknown parameters,'' {\em
  Automatica}, vol.~44, no.~1, pp.~119--127, 2008.

\bibitem{terra2014optimal}
M.~H. Terra, J.~P. Cerri, and J.~Y. Ishihara, ``Optimal robust linear quadratic
  regulator for systems subject to uncertainties,'' {\em IEEE Transactions on
  Automatic Control}, vol.~59, no.~9, pp.~2586--2591, 2014.

\bibitem{qian2020dual}
F.~Qian, X.~Zhang, L.~Liu, and G.~Xie, ``Dual control for stochastic linear
  mimo systems with parameter uncertainty,'' {\em IEEE Access}, vol.~8,
  pp.~41860--41869, 2020.

\bibitem{falkinger2008limited}
J.~Falkinger, ``Limited attention as a scarce resource in information-rich
  economies,'' {\em The Economic Journal}, vol.~118, no.~532, pp.~1596--1620,
  2008.

\bibitem{gargano2018does}
A.~Gargano and A.~G. Rossi, ``Does it pay to pay attention?,'' {\em The Review
  of Financial Studies}, vol.~31, no.~12, pp.~4595--4649, 2018.

\bibitem{peng2006investor}
L.~Peng and W.~Xiong, ``Investor attention, overconfidence and category
  learning,'' {\em Journal of Financial Economics}, vol.~80, no.~3,
  pp.~563--602, 2006.

\bibitem{andrei2020dynamic}
D.~Andrei and M.~Hasler, ``Dynamic attention behavior under return
  predictability,'' {\em Management Science}, vol.~66, no.~7, pp.~2906--2928,
  2020.

\bibitem{zhang2022mean}
Y.~Zhang, Z.~Jin, J.~Wei, and G.~Yin, ``Mean--variance portfolio selection with
  dynamic attention behavior in a hidden markov model,'' {\em Automatica},
  vol.~146, p.~110629, 2022.

\bibitem{wang2022information}
Y.~Wang, D.~Wang, and C.~Hou, ``Information acquisition and asset allocation
  with unknown income growth,'' {\em Economics Letters}, vol.~213, p.~110364,
  2022.

\bibitem{fama1993common}
E.~F. Fama and K.~R. French, ``Common risk factors in the returns on bonds and
  stocks,'' {\em Journal of Financial Economics}, vol.~33, no.~1, pp.~3--53,
  1993.

\bibitem{herskovic2016common}
B.~Herskovic, B.~Kelly, H.~Lustig, and S.~Van~Nieuwerburgh, ``The common factor
  in idiosyncratic volatility: Quantitative asset pricing implications,'' {\em
  Journal of Financial Economics}, vol.~119, no.~2, pp.~249--283, 2016.

\bibitem{garleanu2013dynamic}
N.~G{\^a}rleanu and L.~H. Pedersen, ``Dynamic trading with predictable returns
  and transaction costs,'' {\em The Journal of Finance}, vol.~68, no.~6,
  pp.~2309--2340, 2013.

\bibitem{Arnott2016}
R.~D. Arnott, N.~Beck, and V.~Kalesnik, ``Timing 'smart beta' strategies? of
  course! buy low, sell high!,'' {\em Available at SSRN 3040956}, 2016.

\bibitem{Haddad}
V.~Haddad, S.~Kozak, and S.~Santosh, ``Factor timing,'' {\em The Review of
  Financial Studies}, vol.~33, no.~5, pp.~1980--2018, 2020.

\bibitem{arnott2023factor}
R.~D. Arnott, V.~Kalesnik, and J.~T. Linnainmaa, ``Factor momentum,'' {\em The
  Review of Financial Studies}, vol.~36, no.~8, pp.~3034--3070, 2023.

\bibitem{Favero2021macro}
C.~A. Favero, A.~Melone, and A.~Tamoni, ``Macro trends and factor timing,''
  {\em Available at SSRN 3940452}, 2021.

\bibitem{Bertsekas}
D.~Bertsekas, {\em Dynamic programming and optimal control: Volume I}, vol.~1.
\newblock Athena scientific, 2012.

\end{thebibliography}

\end{document}